\documentclass[psamsfonts, 12pt, a4, leqno]{amsart}

\usepackage{amssymb, latexsym}

\newtheorem{Theorem}{Theorem}

\newtheorem{Lemma}[Theorem]{Lemma}
\newtheorem{Proposition}[Theorem]{Proposition}

\usepackage
{hyperref}

\begin{document}

\title[Discrete mean square of the Riemann zeta-function] {Discrete mean square of the Riemann zeta-function over imaginary parts of its zeros}

\author{ Ram\={u}nas Garunk\v{s}tis}
\address{Ram\={u}nas Garunk\v{s}tis \\
Department of Mathematics and Informatics, Vilnius University \\
Naugarduko 24, 03225 Vilnius, Lithuania}
\email{ramunas.garunkstis@mif.vu.lt}
\urladdr{www.mif.vu.lt/~garunkstis}
\thanks{ The first author is supported
	by a grant No. MIP-049/2014
	from the Research Council of Lithuania.}

\author{Antanas Laurin\v cikas}
\address{Antanas Laurin\v cikas \\
	Department of Mathematics and Informatics, Vilnius University \\
Naugarduko 24, 03225 Vilnius, Lithuania}
\email{antanas.laurincikas@mif.vu.lt}
\urladdr{}

\setcounter{equation}{0}

\date{}

\begin{abstract}
 Assume the Riemann hypothesis. On the right-hand side of the critical strip, we obtain an asymptotic  formula for the discrete mean square of the Riemann zeta-function over imaginary parts of its zeros.
\end{abstract}

\maketitle

\section{Introduction}

Let $s=\sigma+it$ be a complex variable. In this paper, $T$ always tends to plus infinity.

Let $N(T)$ denote the number of zeros of the Riemann zeta-function $\zeta(s)$ in the region $0\le\sigma\le1$, $0<t\le T$. The Riemann-von Mangoldt formula states (Tichmarsh \cite[Theorem 9.4]{titch}) that
\begin{align}\label{vonmang}
N(T)=\frac{T}{2\pi}\log\frac{T}{2\pi e}+O(\log T).
\end{align}
Let $\rho=\beta+i\gamma$ denote a non-real zero of $\zeta(s)$. The Riemann hypothesis (RH) states that $\beta=1/2$ for all non-real zeros of the Riemann zeta-function. We prove the following two theorems.
\begin{Theorem}\label{prop1}
Assume RH. Let $\sigma>1/2$. Let  $A>0$ be as large as we like. Then there is $\eta=\eta(A)>0$ such that, for $|t|\le T^A$,
	\begin{align*}
	\sum_{0<\gamma\le T}|\zeta(s+i\gamma)|^2=\zeta(2\sigma)\frac{T}{2\pi}\log\frac{T}{2\pi e}+\zeta(2\sigma)\Re\left(\frac{\zeta'}{\zeta}(s+1/2)\right)\frac{T}{\pi} +O_{\sigma,A}(T^{1-\eta}).
	\end{align*}
\end{Theorem}
\begin{Theorem}\label{th2}
Assume RH. Let $\sigma>1/2$. Then, for any  $\varepsilon>0$,
	\begin{align*}
	\sum_{0<\gamma\le T}|\zeta(s+i\gamma)|^2\ll_{\sigma,\varepsilon} T\log T+|t|^\varepsilon,
	\end{align*}
	uniformly in $t$.
\end{Theorem}
In Garunk\v stis and Laurin\v cikas \cite{gl}, we use Theorem \ref{th2} in order to study the discrete universality of the Riemann zeta-function over the imaginary parts of its zeros. Informally speaking, this means that a wide class of analytic functions can be approximated by shifts $\zeta(s+i\gamma)$. We note that the discrete universality was proposed by Reich \cite{rei}. It was developed by Bagchi \cite{bag},  Sander and Steuding \cite{sast}. These authors investigated the approximation of analytic functions by shifts $\zeta(s+i\tau)$, where $\tau$ takes values from arithmetic progression $\{kh : k=0,1,2,\dots\}$ with $h>0$ fixed. Instead of arithmetic progressions, Dubickas and Laurin\v cikas \cite{dl} considered the set $\{k^\alpha h : k=0,1,2,\dots\}$ with $0<\alpha<1$ fixed.

The related discrete mean square  was considered by Gonek \cite{gon84}.  He proved,  assuming RH, that for real $\alpha$, $|\alpha|\le (1/4\pi)\log (T/2\pi)$,
$$\sum_{0<\gamma\le T}\left|\zeta\left(\frac12+i\left(\gamma+\frac{2\pi\alpha}{\log (T/2\pi)} \right)\right)\right|^2=\left(1-\left(\frac{\sin (\pi\alpha)}{\pi\alpha}\right)^2\right)\frac{T}{2\pi}\log^2 T+O(T\log T).$$
 The error term in the last formula was improved (on RH) by Fujii \cite{Fujii95}, see also Conrey and Snaith \cite[Section 7.3]{cs}, where this formula is investigated using the ratio conjecture.  Ivi\'c \cite{ivi} obtained that
$$
\sum_{0<\gamma\le T}\left|\zeta\left(\frac12+i\gamma\right)\right|^2\ll T\log^2 T\log\log T
$$
unconditionally.

In the next section, we prove Theorems \ref{prop1} and \ref{th2}. In Section \ref{conclrem}, we discuss several discrete mean square results for Dirichlet $L$-functions.

\section{Proofs}
  In the proofs of Theorems \ref{prop1} and \ref{th2},  we will use the approximation of $\zeta(s)$ by a finite sum and the uniform version of Landau's formula (see Lemmas \ref{funcEq} and \ref{A} below).

 \begin{Lemma}\label{funcEq}
 	Assume RH. Let $\sigma>1/2$ and $t>0$. Then, for any given positive number $\delta$, there is $\lambda=\lambda(\delta,\sigma)>0$ such that
	\begin{align*}
 	\zeta(s) =\sum_{n\le t^\delta}
 	\frac{1}{n^{s}}+r(s),
 	\end{align*}
	where  $r(s)\ll t^{-\lambda}$.
 \end{Lemma}
 \proof
 The lemma follows from Tichmarsh \cite[Theorem 13.3]{titch}.
 \endproof

 \begin{Lemma}\label{A}
 	Assume RH. Let   $x,T>1$. Then
	\begin{align*}
\sum_{0<\gamma\leq T}x^{i\gamma } =& -\frac{T}{2\pi}\frac{\Lambda(x)}{\sqrt{x}}+O\left(\sqrt{x}\log(2xT)\log\log (3x) \right)
\\
&+O\left(\frac{\log x}{\sqrt{x}} \min\left(T, \frac{x}{\langle x\rangle}\right)\right)+O\left(\frac{\log (2T)}{\sqrt{x}} \min\left(T, \frac{1}{\log x}\right)\right),
 	\end{align*}
where $\langle x\rangle$ denotes the distance from $x$ to the nearest prime power other than $x$ itself.
 \end{Lemma}
\proof
Under RH, the lemma  follows immediately from Gonek  \cite{gon1} and \cite{gon2}. Note that stronger forms of Landau's formula are obtained by Fujii \cite{Fujii1}, \cite{Fujii2} (under RH), also by Ford and  Zaharescu \cite{fz}.

\endproof

The following lemma will be useful.
\begin{Lemma}\label{gammat}
Let $0<\lambda<1/2$. Let $t>0$ or $-t>T$. Then
\begin{align*}
\sum_{0<\gamma\le T}(\gamma+t)^{-2\lambda}\ll_\lambda T|T+t|^{-2\lambda}\log T\ll_\lambda T|T+t|^{-\lambda}
\end{align*}
uniformly in $t$.
\end{Lemma}
\proof
By partial summation and by the Riemann-von Mangoldt formula \eqref{vonmang}, we have
\begin{align*}
&\sum_{0<\gamma\le T}(\gamma+t)^{-2\lambda}=\frac{N(T)}{(T+t)^{2\lambda}}+2\lambda\int_0^T\frac{N(x)dx}{(x+t)^{2\lambda+1}}
\\& \ll\frac{T\log T}{(T+t)^{2\lambda}} +\log T \int_1^T\frac{xdx}{(x+t)^{2\lambda+1}}.
\end{align*}
Let $I=\int_1^Tx(x+t)^{-2\lambda-1}dx$. Then
\begin{align}\label{Iexpr}
I=\frac{1}{1-2\lambda}\left((T+t)^{1-2\lambda}-(1+t)^{1-2\lambda}\right)+\frac{t}{2\lambda}\left((T+t)^{-2\lambda}-(1+t)^{-2\lambda}\right).
\end{align}
If $0\le t\le T$ then
\begin{align*}
I\ll (T+t)^{1-2\lambda}+t(1+t)^{-2\lambda}\ll T(T+t)^{-2\lambda}.
\end{align*}
If $|t|>T$  then, using the formula \eqref{Iexpr} and Taylor series, we get
\begin{align*}
I=&\frac{t^{1-2\lambda}}{1-2\lambda}\left(\left(1+\frac{T}{t}\right)^{1-2\lambda}-\left(1+\frac{1}{t}\right)^{1-2\lambda}\right)
\\&
+\frac{t^{1-2\lambda}}{2\lambda}\left(\left(1+\frac{T}{t}\right)^{-2\lambda}-\left(1+\frac{1}{t}\right)^{-2\lambda}\right)
\\=&
\frac{t^{1-2\lambda}}{1-2\lambda}\left((1-2\lambda)\frac Tt-(1-2\lambda)\frac 1t+O\left(\frac {T^2}{t^2}\right)\right)
\\&+
\frac{t^{1-2\lambda}}{2\lambda}\left(-2\lambda\frac Tt+2\lambda\frac 1t+O\left(\frac {T^2}{t^2}\right)\right)
\\ \ll&
T^2|t|^{-1-2\lambda}\ll T|t|^{-2\lambda}\ll T|T+t|^{-2\lambda}.
\end{align*}
Lemma \ref{gammat} is proved.
\endproof

To prove Theorems \ref{prop1} and \ref{th2}, we consider two cases: $t\ge 0$ and $t < 0$. These cases correspond to Propositions \ref{Ttpositive} and \ref{0<-t<T} below.

\begin{Proposition}\label{Ttpositive}
Assume RH. Let $\sigma>1/2$, $t\ge0$, and $0<\delta<1$. Then there is a positive number $\lambda=\lambda(\delta,\sigma)$ such that
	\begin{align*}
	&\sum_{0<\gamma\le T}|\zeta(s+i\gamma)|^2=\zeta(2\sigma)\frac{T}{2\pi}\log\frac{T}{2\pi e}+\zeta(2\sigma)\Re\left(\frac{\zeta'}{\zeta}(s+1/2)\right)\frac{T}{\pi}
\\&+
O\left(T(T+t)^{-\delta(\sigma-1/2)}+T(T+t)^{-\lambda}+(T+t)^{2\delta}+T^{1/2}(T+t)^{\delta-\lambda}\right).
	\end{align*}
\end{Proposition}

 \proof
In view of Lemma \ref{funcEq}, we have
\begin{align}\label{AR}
	&\sum_{0<\gamma\le T}|\zeta(\sigma+i\gamma+it)|^2
	=
	\sum_{\gamma\le T}\sum_{n,m\le (\gamma+t)^\delta}\frac{(m/n)^{i\gamma+it}}{(mn)^\sigma}
\\&
+\sum_{0<\gamma\le T}\left(|r(\sigma+i\gamma+it)|^2+2\Re\sum_{n\le(\gamma+t)^\delta}\frac{\overline{r(\sigma+i\gamma+it)}}{n^{\sigma+i\gamma+it}}\right)
\nonumber
	\\&=:
	A+R.\nonumber
\end{align}
By Lemmas \ref{funcEq}, \ref{gammat},
  and by the Cauchy-Schwarz inequality, there exists a number $0<\lambda<1/2$ such that
\begin{align}\label{R}
R&\ll \sum_{0<\gamma\le T}|r(\sigma+i\gamma+it)|^2 +\left(A\sum_{0<\gamma\le T}|r(\sigma+i\gamma+it)|^2\right)^{1/2}
\\&\ll
 T(T+t)^{-\lambda}+A^{1/2}T^{1/2} (T+t)^{-\lambda}\log^{1/2} T.\nonumber
\end{align}
We rewrite the term $A$ in the following way
\begin{align}\label{multiplesum}
	A&=
	\sum_{\gamma\le T}\sum_{n\le(\gamma+t)^\delta}\frac{1}{n^{2\sigma}}
+
2\Re\sum_{\gamma\le T}\sum_{n<m\le (\gamma+t)^\delta}\frac{(m/n)^{i\gamma+it}}{(mn)^\sigma}
	\\&=:
	A_1+2\Re A_2.\nonumber
\end{align}
By the Riemann-von Mangoldt formula \eqref{vonmang}, we see that
\begin{align}\label{FirstDoublesummOfA}
A_1=\zeta(2\sigma)\frac{T}{2\pi}\log\frac{T}{2\pi e}+O\left(T(T+t)^{-\delta(2\sigma-1)}\log(T+t)+\log T\right).
\end{align}

The sum $A_2$ requires longer consideration. Changing the order of the summation, we obtain
\begin{align}\label{Epr5}
A_2=\sum_{\gamma\le T}\sum_{n<m\le (\gamma+t)^\delta}\frac{(m/n)^{i\gamma+it}}{(mn)^\sigma}=\sum_{n<m\le (T+t)^\delta}\sum_{\max(0,m^{1/\delta}-t)<\gamma\le T }\frac{(m/n)^{i\gamma+it}}{(mn)^\sigma}.
\end{align}
On the right-hand side of the last equality, we will use Lemma \ref{A} for the inner sum.  We have
\begin{align*}
\left\langle\frac{m}{n}\right\rangle\ge\min_{{d\ne m/n\atop d\in\mathbb N}}\left|\frac{m}{n}-d\right|=\min_{{d\ne m/n\atop d\in\mathbb N}}\frac1n\left|m-dn\right|\ge\frac1n
\end{align*}
and
\begin{align*}
\log \frac mn=
\log
\left(
1+\frac{m-n}n
\right)
\gg \frac1{n}.
\end{align*}
Then by Lemma \ref{A}, for $n<m\le(T+t)^{\delta}$, we get
\begin{align}\label{m/n}
	\sum_{\max(0,m^{1/\delta}-t)<\gamma\le T }(m/n)^{i\gamma}
=&-\frac{T-\max(0,m^{1/\delta}-t)}{2\pi}\frac{\Lambda(m/n)}{\sqrt{m/n}}
	\\&
	+O\left( (T+t)^{\delta}\log(T+t)\right).\nonumber
		\end{align}
Further,
\begin{align}\label{lambdamn}
&\sum_{n<m\le(T+t)^\delta}\frac{(m/n)^{it}\Lambda(m/n)}{(mn)^\sigma\sqrt{m/n}}
=
\sum_{n<p^kn\le(T+t)^\delta}\frac{p^{kit}\Lambda(p^k)}{(p^kn^2)^\sigma\sqrt{p^k}}
\\&=
\sum_{n\le(T+t)^\delta}\frac{1}{n^{2\sigma}}\sum_{j\le(T+t)^\delta/n}\frac{j^{it}\Lambda(j)}{j^{\sigma+1/2}}\nonumber
\\&=
-\sum_{n\le(T+t)^\delta}\frac{1}{n^{2\sigma}}\left(\frac{\zeta'}{\zeta}(\sigma-it+1/2)+\sum_{j>(T+t)^\delta/n}\frac{j^{it}\Lambda(j)}{j^{\sigma+1/2}}\right)\nonumber
\\&=
-\sum_{n\le(T+t)^\delta}\frac{1}{n^{2\sigma}}\left(\frac{\zeta'}{\zeta}(\sigma-it+1/2)+O\left(\left(\frac{(T+t)^\delta}{n}\right)^{1/2-\sigma}\right)\right).\nonumber
\end{align}
In view of
\begin{align*}
\sum_{n\le(T+t)^\delta}\frac{1}{n^{\sigma+1/2}}\ll1\quad\text{and}\quad \sum_{n>(T+t)^\delta}\frac{1}{n^{2\sigma}}\ll (T+t)^{\delta(1-2\sigma)},
\end{align*}
we get
\begin{align}\label{MainTerm}
&\sum_{n<m\le(T+t)^\delta}\frac{(m/n)^{it}\Lambda(m/n)}{(mn)^\sigma\sqrt{m/n}}
\\&=
-\zeta(2\sigma)\frac{\zeta'}{\zeta}(\sigma-it+1/2)+O((T+t)^{-\delta(\sigma-1/2)}).\nonumber
\end{align}
We continue to consider the right-hand side of the formula \eqref{m/n}. Reasoning similarly as in (\ref{lambdamn}), we obtain
\begin{align}\label{S}
S&:=\left|\sum_{n<m\le(T+t)^\delta}\frac{(m/n)^{it}\max(0,m^{1/\delta}-t)\Lambda(m/n)}{(mn)^\sigma\sqrt{m/n}}\right|
\\&=
\left|\sum_{n<p^kn\le(T+t)^\delta}\frac{p^{kit}\max(0,(p^kn)^{1/\delta}-t)\Lambda(p^k)}{(p^kn^2)^\sigma\sqrt{p^k}}\right|\nonumber
\\&=
\left|\sum_{n\le(T+t)^\delta}\frac{1}{n^{2\sigma}}\sum_{t^\delta/n<j\le(T+t)^\delta/n}\frac{j^{it}((jn)^{1/\delta}-t)\Lambda(j)}{j^{\sigma+1/2}}\right|\nonumber
\\&\le
\sum_{n\le(T+t)^\delta}n^{1/\delta-2\sigma}\sum_{t^\delta/n<j\le(T+t)^\delta/n}j^{1/\delta-\sigma-1/2}\Lambda(j).\nonumber
\end{align}
In light of
\begin{align*}
\sum_{t^\delta/n<j\le(T+t)^\delta/n}j^{1/\delta-\sigma-1/2}\Lambda(j)
\ll&
 \frac{(T+t)^{1-\delta(\sigma-1/2)}-t^{1-\delta(\sigma-1/2)}}{n^{1/\delta-\sigma+1/2}}
 \\ <&
 \frac{T(T+t)^{-\delta(\sigma-1/2)}}{n^{1/\delta-\sigma+1/2}},
\end{align*}
we see that
\begin{align}\label{FinalS}
S\ll T(T+t)^{-\delta(\sigma-1/2)}.
\end{align}
Hence,  putting together formulas \eqref{Epr5}--\eqref{FinalS},  we get
\begin{align*}
2\Re A_2
&=
\zeta(2\sigma)\Re\left(\frac{\zeta'}{\zeta}(s+1/2)\right)\frac{T}{\pi}
+
O\left(T(T+t)^{-\delta(\sigma-1/2)}\right)\nonumber
\\&
\phantom{=}+
O\left(((T+t)^{2\delta(1-\sigma)} +1)(T+t)^{\delta}\log(T+t)\right).
\end{align*}

By the last formula together with formulas \eqref{multiplesum} and \eqref{FirstDoublesummOfA}, we have
\begin{align}\label{finalA}
A=&\zeta(2\sigma)\frac{T}{2\pi}\log\frac{T}{2\pi e}+\zeta(2\sigma)\Re\left(\frac{\zeta'}{\zeta}(s+1/2)\right)\frac{T}{\pi}
\\&+
O\left(T(T+t)^{-\delta(\sigma-1/2)}+(T+t)^{2\delta}\right).\nonumber
\end{align}
Then  by \eqref{R}, we see that $R\ll T(T+t)^{-\lambda/2}+T^{1/2}(T+t)^{\delta-\lambda/2}$.
From this and formulas \eqref{AR}, \eqref{finalA},  replacing $\lambda/2$ by $\lambda$, we obtain    Proposition \ref{Ttpositive}.

 \endproof

\begin{Proposition}\label{0<-t<T}
Assume RH. Let $\sigma>1/2$ and  let
\begin{align*}
	\sum_{\gamma\le T}|\zeta(s+i\gamma)|^2
=
\zeta(2\sigma)\frac{T}{2\pi}\log\frac{T}{2\pi e}+\zeta(2\sigma)\Re\left(\frac{\zeta'}{\zeta}(s+1/2)\right)\frac{T}{\pi}
+E.
	\end{align*}
If  $0<-t\le 2T$, then there is a positive number $\eta=\eta(\sigma)$ such that
$$E\ll T^{1-\eta}.$$

Moreover, let $0<\delta\le1/2$.  If $-t>2T>0$, then there is a positive number $\lambda=\lambda(\delta,\sigma)$ such that
$$E\ll T(-t)^{-\delta(\sigma-1/2)}+T(-t)^{-\lambda}+(-t)^{2\delta}
+
T^{1/2}(-t)^{\delta-\lambda}.$$
\end{Proposition}

\proof
This proof differs from the proof of Proposition \ref{Ttpositive}  in the way that now the number $\gamma+t$ could be negative and its absolute value could be small. 

Let $\varepsilon>0$. First we consider the case  $0<-t\le T+T^\varepsilon$. If $|\gamma+t|$ is small, then, using the bound $\zeta(s)\ll t^{\varepsilon/3}$ (Titchmarsh \cite[formula (14.2.5)]{titch}) and the Riemann-von Mangoldt formula \eqref{vonmang}, we have
\begin{align*}
	\sum_{{0<\gamma\le T\atop |t+\gamma|\le T^\varepsilon}}|\zeta(s+i\gamma)|^2\ll  T^{2\varepsilon}.
	\end{align*}
Moreover, we separate the terms with negative $\gamma+t$  into a different sum. Then
\begin{align}\label{sumsplit}
\sum_{0<\gamma\le T}|\zeta(\sigma+i\gamma+it)|^2
=&
\sum_{0<\gamma\le -t-T^\varepsilon}|\zeta(\sigma+i\gamma+it)|^2
\\&+
\sum_{T^\varepsilon -t<\gamma\le T}|\zeta(\sigma+i\gamma+it)|^2
+
O(T^{2\varepsilon}).\nonumber
\end{align}
Here and later we assume that the empty sum is equal to zero.

We consider the first sum on the right-hand side of the last equality. Lemma~\ref{gammat} gives that
\begin{align}\label{nexttoAR2}
\sum_{0<\gamma\le -t-T^\varepsilon}(\gamma+t)^{-2\lambda}\ll_\lambda T^{1-2\varepsilon\lambda}\log T\ll_\lambda T^{1-\varepsilon\lambda}.
\end{align}
Then by Lemma \ref{funcEq}, reasoning similarly as in formulas \eqref{AR}--\eqref{Epr5}, we have that there exists a number $0<\lambda<1/2$ such that
\begin{align}\label{AR2}
&\sum_{0<\gamma\le -t-T^\varepsilon}|\zeta(\sigma+i\gamma+it)|^2
=:
\sum_{\gamma\le -t-T^\varepsilon}\sum_{n\le(-t-\gamma)^\delta}\frac{1}{n^{2\sigma}}
\\&+
2\Re\left(\sum_{n<m\le (-t)^\delta}\sum_{\gamma\le \min(-T^\varepsilon, -m^{1/\delta})-t}\frac{(m/n)^{i\gamma+it}}{(mn)^\sigma}\right)+R
=:
A+R,\nonumber
\end{align}
with
\begin{align}\label{rll}
R\ll A^{1/2} T^{1/2-\varepsilon\lambda}\log^{1/2} T+T^{1-\varepsilon\lambda}.
\end{align}
By the Riemann-von Mangoldt formula \eqref{vonmang} and  partial summation
\begin{align}\label{A=}
A=&\zeta(2\sigma)\frac{-t}{2\pi}\log\frac{-t}{2\pi e}+2\Re\left(\sum_{n<m\le (-t)^\delta}\sum_{\gamma\le \min(-T^\varepsilon, -m^{1/\delta})-t}\frac{(m/n)^{i\gamma+it}}{(mn)^\sigma}\right)
\\&+
O\left((T^{1-\varepsilon\delta(2\sigma-1)}+T^\varepsilon)\log T\right).\nonumber
\end{align}
 We consider the double sum of the last formula. Using Lemma \ref{A}, for $n<m\le(-t)^{\delta}$, we find
\begin{align}\label{minteps}
	&\sum_{\gamma\le\min(-T^\varepsilon, -m^{1/\delta})-t}(m/n)^{i\gamma}=-\frac{\min(-T^\varepsilon, -m^{1/\delta})-t}{2\pi}\frac{\Lambda(m/n)}{\sqrt{m/n}}
\\&
+O\left( T^{\delta}\log T\right).\nonumber
		\end{align}
Reasoning similarly as in formulas \eqref{lambdamn} and \eqref{MainTerm}, we get
\begin{align}\label{asinformula}
&t\sum_{n<m\le(-t)^\delta}\frac{(m/n)^{it}\Lambda(m/n)}{(mn)^\sigma\sqrt{m/n}}
=
t\sum_{n\le(-t)^\delta}\frac{1}{n^{2\sigma}}\sum_{j\le(-t)^\delta/n}\frac{j^{it}\Lambda(j)}{j^{\sigma+1/2}}
\\&=
-t\zeta(2\sigma)\frac{\zeta'}{\zeta}(\sigma-it+1/2)+O((-t)^{1-\delta(\sigma-1/2)}).\nonumber
\\&=
-t\zeta(2\sigma)\frac{\zeta'}{\zeta}(\sigma-it+1/2)+O(T^{1-\delta(\sigma-1/2)}+1).\nonumber
\end{align}
Further,
\begin{align}\label{followingformula}
\sum_{n<m\le(-t)^\delta}\frac{m^{1/\delta}\Lambda(m/n)}{(mn)^\sigma\sqrt{m/n}}
&=
\sum_{n\le(-t)^\delta}n^{1/\delta-2\sigma}\sum_{j\le(-t)^\delta/n}j^{1/\delta-\sigma-1/2}\Lambda(j)
\\&\ll
(-t)^{1-\delta(\sigma-1/2)}+1\ll
T^{1-\delta(\sigma-1/2)}+1.\nonumber
\end{align}
By the inequality $\min(-T^\varepsilon, -m^{1/\delta})\le T^\varepsilon+m^{1/\delta}$, formulas \eqref{minteps}, \eqref{asinformula}, and \eqref{followingformula} we obtain
\begin{align}\label{2Resum}
&2\Re\left(\sum_{n<m\le (-t)^\delta}\sum_{\gamma\le \min(-T^\varepsilon, -m^{1/\delta})-t}\frac{(m/n)^{i\gamma+it}}{(mn)^\sigma}\right)
\\&	=
\zeta(2\sigma)\Re\left(\frac{\zeta'}{\zeta}(s+1/2)\right)\frac{-t}{\pi}
+
O(T^{1-\delta(\sigma-1/2)}+ T^{2\delta}+T^\varepsilon).\nonumber
\end{align}
This and the equality \eqref{A=} yield that
\begin{align}\label{lasta}
A=&\zeta(2\sigma)\frac{-t}{2\pi}\log\frac{-t}{2\pi e}+\zeta(2\sigma)\Re\left(\frac{\zeta'}{\zeta}(s+1/2)\right)\frac{-t}{\pi}
\\&
+O(T^{1-\varepsilon\delta(\sigma-1/2)}+ T^{2\delta}+T^\varepsilon\log T).\nonumber
\end{align}
Then in view of \eqref{rll} we see that $R\ll T^{1-\varepsilon\lambda/2}$. By the last bound together with formulas \eqref{AR2} and \eqref{lasta}
 we get 
\begin{align}\label{firstsum}
\sum_{0<\gamma\le -t-T^\varepsilon}&|\zeta(\sigma+i\gamma+it)|^2
=\zeta(2\sigma)\frac{-t}{2\pi}\log\frac{-t}{2\pi e}+\zeta(2\sigma)\Re\left(\frac{\zeta'}{\zeta}(s+1/2)\right)\frac{-t}{\pi}\nonumber
\\&+
O\left(T^{1-\varepsilon\lambda/2}+T^{1-\varepsilon\delta(\sigma-1/2)}+T^{2\delta}+T^\varepsilon\log T\right),
\end{align}
where $0<-t\le T+T^\varepsilon$.

We turn to the next 
sum on the right-hand side of the formula \eqref{sumsplit}. That is, we consider the sum
\begin{align}\label{nextinsumsplit}
\sum_{T^\varepsilon -t<\gamma\le T}|\zeta(\sigma+i\gamma+it)|^2.
\end{align}
Note that the last sum is empty if $T-T^\varepsilon<-t\le T+T^\varepsilon$. We therefore assume that
\begin{align}\label{wetherefore}
0<-t\le T-T^\varepsilon
\end{align}
in the formula \eqref{nextinsumsplit}.
By the Riemann-von Mangoldt formula \eqref{vonmang},
\begin{align*}
\sum_{T^\varepsilon-t<\gamma\le T}(\gamma+t)^{-2\lambda}\ll_\lambda T^{1-2\varepsilon\lambda}\log T\ll_\lambda T^{1-\varepsilon\lambda}.
\end{align*}
Then
\begin{align}\label{A1R1}
&\sum_{T^\varepsilon-t<\gamma\le T}|\zeta(\sigma+i\gamma+it)|^2=:\sum_{T^\varepsilon-t<\gamma\le T}\sum_{n\le(t+\gamma)^\delta}\frac{1}{n^{2\sigma}}
\\&
+2\Re\left(\sum_{n<m\le (T+t)^\delta}\sum_{\max(T^\varepsilon, m^{1/\delta})-t<\gamma\le T }\frac{(m/n)^{i\gamma+it}}{(mn)^\sigma}\right)+R'\nonumber
\\&=:
	A'+R', \nonumber
\end{align}
where
$$R'\ll A'^{1/2} T^{1/2-\varepsilon\lambda}\log^{1/2} T+T^{1-\varepsilon\lambda}.$$
Further,
\begin{align}\label{T-|t|}
&A'=\zeta(2\sigma)\left(\frac{T}{2\pi}\log\frac{T}{2\pi e}+\frac{t}{2\pi}\log\frac{-t}{2\pi e}\right)
\\&
+2\Re\left(\sum_{n<m\le (T+t)^\delta}\sum_{\max(T^\varepsilon, m^{1/\delta})-t<\gamma\le T }\frac{(m/n)^{i\gamma+it}}{(mn)^\sigma}\right)\nonumber
\\&
+O\left((T^{1-\varepsilon\delta(2\sigma-1)}+T^\varepsilon)\log T\right).\nonumber
\end{align}
 In view of inequalities \eqref{wetherefore} we have $T+t\ge T^\varepsilon$. We split the double sum in the equation \eqref{T-|t|} in the following way.
\begin{align}\label{onemore}
&\sum_{n<m\le (T+t)^\delta}\sum_{\max(T^\varepsilon, m^{1/\delta})-t<\gamma\le T}\frac{(m/n)^{i\gamma+it}}{(mn)^\sigma}
=\sum_{n<m\le  T^{\varepsilon\delta}}\sum_{T^\varepsilon-t<\gamma\le T}\frac{(m/n)^{i\gamma+it}}{(mn)^\sigma}
\\&
+\sum_{T^{\varepsilon\delta}<m\le (T+t)^\delta}\sum_{n<m}\sum_{ m^{1/\delta}-t<\gamma\le T}\frac{(m/n)^{i\gamma+it}}{(mn)^\sigma}=:C'+D'\nonumber
\end{align}
For $n<m\le(T+t)^{\delta}$, Lemma \ref{A} yields
\begin{align}\label{lema4yields}
	\sum_{\max(T^\varepsilon, m^{1/\delta})-t<\gamma\le T}(m/n)^{i\gamma}=-\frac{T+t-\max(T^\varepsilon, m^{1/\delta})}{2\pi}\frac{\Lambda(m/n)}{\sqrt{m/n}}
	+O\left( T^{\delta}\log T\right).
\end{align}
Reasoning similarly as in formulas \eqref{lambdamn} and \eqref{MainTerm}, we get
\begin{align}\label{termC'}
C'=&-\frac{T+t-T^\varepsilon}{2\pi}\sum_{n<m\le T^{\varepsilon\delta}}\frac{(m/n)^{it}\Lambda(m/n)}{(mn)^\sigma\sqrt{m/n}}+O\left(T^{2\delta}\right)
\\=&
\frac{T+t}{2\pi}\zeta(2\sigma)\frac{\zeta'}{\zeta}(\sigma-it+1/2)+O\left(T^{1-\varepsilon\delta(\sigma-1/2)}+T^{2\delta}+T^\varepsilon\right)\nonumber
\end{align}

We turn to the term $D'$ defined by the formula \eqref{onemore}.  In view of \eqref{lema4yields} we obtain 
\begin{align}\label{d'=}
D'=-\sum_{T^{\varepsilon\delta}<m\le (T+t)^\delta}\sum_{n<m}\frac{T+t- m^{1/\delta}}{2\pi}\frac{(m/n)^{it}\Lambda(m/n)}{(mn)^\sigma\sqrt{m/n}}
	+O\left( T^{2\delta}\right).
\end{align}
 In formula \eqref{d'=} we have that $T+t-m^{1/\delta}<T$. Thus
\begin{align}\label{d'1d'2}
&D'\ll \sum_{T^{\varepsilon\delta}<m\le (T+t)^\delta}\sum_{n<m}\frac{T\Lambda(m/n)}{(mn)^\sigma\sqrt{m/n}}+O\left( T^{2\delta}\right)
\\&
=\sum_{n<(T+t)^\delta}\sum_{\max(n,T^{\varepsilon\delta})<m\le(T+t)^\delta}\frac{T\Lambda(m/n)}{(mn)^\sigma\sqrt{m/n}}+O\left( T^{2\delta}\right)\nonumber
\\&
=\sum_{n\le T^{\varepsilon\delta}}\sum_{T^{\varepsilon\delta}<m\le (T+t)^\delta}\frac{T\Lambda(m/n)}{(mn)^\sigma\sqrt{m/n}}\nonumber
\\&+\sum_{T^{\varepsilon\delta}<n<(T+t)^\delta}\sum_{n<m\le (T+t)^\delta}\frac{T\Lambda(m/n)}{(mn)^\sigma\sqrt{m/n}}+O\left( T^{2\delta}\right)\nonumber
\\&
=:D'_1+D'_2+O\left( T^{2\delta}\right).\nonumber
\end{align}
Reasoning  as in formulas \eqref{lambdamn} and \eqref{MainTerm}, we get
\begin{align}\label{d'1d'2ll}
D'_1, D'_2\ll T^{1-\varepsilon\delta(\sigma-1/2)}.
\end{align}
Formulas \eqref{d'1d'2} and \eqref{d'1d'2ll} yield the bound
\begin{align}\label{yieldboundd'}
D'\ll T^{1-\varepsilon\delta(\sigma-1/2)}+T^{2\delta}.
\end{align}

Formulas \eqref{onemore}, \eqref{termC'}, and \eqref{yieldboundd'} give
\begin{align}\label{T+t1}
&2\Re\left(\sum_{n<m\le (T+t)^\delta}\sum_{\max(T^\varepsilon, m^{1/\delta})-t<\gamma\le T}\frac{(m/n)^{i\gamma+it}}{(mn)^\sigma}\right)
\\&=
\zeta(2\sigma)\Re\left(\frac{\zeta'}{\zeta}(s+1/2)\right)\frac{T+t}{\pi}
+
O\left(T^{1-\delta(\sigma-1/2)}+T^{2\delta}+T^\varepsilon\right).\nonumber
\end{align}
Note that the error terms in formulas \eqref{A1R1},  \eqref{T-|t|}, and \eqref{T+t1} are the same as in corresponding formulas \eqref{AR2}, \eqref{A=}, and \eqref{2Resum}. Therefore, similarly to the derivation of the formula \eqref{firstsum}, using \eqref{A1R1},  \eqref{T-|t|}, and \eqref{T+t1}, we get
\begin{align}\label{secondsum}
&\sum_{T^\varepsilon-t<\gamma\le T}|\zeta(\sigma+i\gamma+it)|^2=\zeta(2\sigma)\left(\frac{T}{2\pi}\log\frac{T}{2\pi e}+\frac{t}{2\pi}\log\frac{-t}{2\pi e}\right)
\\&+
\zeta(2\sigma)\Re\left(\frac{\zeta'}{\zeta}(s+1/2)\right)\frac{T+t}{\pi}
\nonumber
\\&+
O\left(T^{1-\varepsilon\lambda/2}+T^{1-\varepsilon\delta(\sigma-1/2)}+T^{2\delta}+T^\varepsilon\log T\right),\nonumber
\end{align}
where $0<-t\le T-T^\varepsilon$. Note again that the sum \eqref{secondsum} is empty for $T-T^\varepsilon<-t\le T+T^\varepsilon$.

Next we consider the case $-t>T+T^\varepsilon$. Lemma \ref{gammat} gives that
\begin{align*}
&\sum_{0<\gamma\le T}(\gamma+t)^{-2\lambda}\\&
\ll_\lambda
\begin{cases}
T^{1-2\varepsilon\lambda}\log T
\ll_\lambda T^{1-\varepsilon\lambda} \quad&\text{if}\quad T+T^\varepsilon<-t\le2T,
\\
T(-t)^{-2\lambda}\log T\ll_\lambda T(-t)^{-\lambda} \quad&\text{if}\quad -t>2T.
\end{cases}
\end{align*}
Therefore
\begin{align}\label{AProp7}
&\sum_{0<\gamma\le T}|\zeta(\sigma+i\gamma+it)|^2=:\sum_{\gamma\le T}\sum_{n\le(-t-\gamma)^\delta}\frac{1}{n^{2\sigma}}
\\&+
2\Re\left(\sum_{n<m\le (-t)^\delta}\sum_{\gamma\le \min(T,-t-m^{1/\delta})}\frac{(m/n)^{i\gamma+it}}{(mn)^\sigma}\right)
+R''
=:A''+R'',\nonumber
\end{align}
where
\begin{align}\label{r''}
R''\ll
\begin{cases}
 A''^{1/2}T^{1/2-\varepsilon\lambda} \log^{1/2} T+T^{1-\varepsilon\lambda} \quad&\text{if}\quad T+T^\varepsilon<-t\le2T,
\\
A''^{1/2}T^{1/2}(-t)^{-\lambda} \log^{1/2} T+T(-t)^{-\lambda} \quad&\text{if}\quad -t>2T.
\end{cases}
\end{align}
Futher
\begin{align}\label{BProp7}
A''=&\zeta(2\sigma)\frac{T}{2\pi}\log\frac{T}{2\pi e}+2\Re\sum_{n<m\le (-t)^\delta}\sum_{\gamma\le \min(T,-t-m^{1/\delta}) }\frac{(m/n)^{i\gamma+it}}{(mn)^\sigma}
\\&+
R_{A''},\nonumber
\end{align}
where 
\begin{align}\label{rall}
R_{A''}\ll
\begin{cases}
 T^{1-\varepsilon\delta(2\sigma-1)}\log T+\log T \quad&\text{if}\quad T+T^\varepsilon<-t\le2T,
\\
T(-t)^{-\delta(2\sigma-1)}\log T+\log T \quad&\text{if}\quad -t>2T.
\end{cases}
\end{align}
We split the double sum in the equation \eqref{BProp7} in the following way.
\begin{align}\label{sumsplit1}
&\sum_{n<m\le (-t)^\delta}\sum_{\gamma\le \min(T,-t-m^{1/\delta}) }\frac{(m/n)^{i\gamma+it}}{(mn)^\sigma}
\\=&
\sum_{n<m\le (-t-T)^\delta}\sum_{\gamma\le T }\frac{(m/n)^{i\gamma+it}}{(mn)^\sigma}\nonumber
\\&+
\sum_{(-t-T)^\delta<m\le (-t)^\delta}\sum_{n<m}\sum_{\gamma\le -t-n^{1/\delta} }\frac{(m/n)^{i\gamma+it}}{(mn)^\sigma}.\nonumber
\\=&:C''+D''.\nonumber
\end{align}
Then 
\begin{align}\label{acd}
A''=\zeta(2\sigma)\frac{T}{2\pi}\log\frac{T}{2\pi e}+2\Re(C''+D'')+R_{A''}.
\end{align}
Next we consider terms $C''$ and $D''$. By Lemma \ref{A}, for $n<m\le(-t)^{\delta}$, we obtain
\begin{align}\label{minT}
	\sum_{\gamma\le \min(T,-t-m^{1/\delta})}(m/n)^{i\gamma}=&-\frac{\min(T,-t-m^{1/\delta})}{2\pi}\frac{\Lambda(m/n)}{\sqrt{m/n}}
\\&
	+O\left( (-t)^{\delta}\log (-t)\right).\nonumber
		\end{align}
By the last formula we get
\begin{align}\label{termC}
&C''=-\frac{T}{2\pi}\sum_{n<m\le (-t-T)^\delta}\frac{(m/n)^{it}\Lambda(m/n)}{(mn)^\sigma\sqrt{m/n}}+R_{C''},
\end{align}
where 
\begin{align*}
R_{C''}\ll
\begin{cases}
 T^{2\delta} \quad&\text{if}\quad T+T^\varepsilon<-t\le2T,
\\
(-t)^{2\delta} \quad&\text{if}\quad -t>2T.
\end{cases}
\end{align*}
 Following the reasoning used in \eqref{lambdamn} and \eqref{MainTerm} we see that
\begin{align}
\sum_{n<m\le (-t-T)^\delta}\frac{(m/n)^{it}\Lambda(m/n)}{(mn)^\sigma\sqrt{m/n}}
=
-\zeta(2\sigma)\frac{\zeta'}{\zeta}(\sigma+1/2-it)+Q,\nonumber
\end{align}
where 
\begin{align*}
Q\ll
\begin{cases}
 T^{-\varepsilon\delta(\sigma-1/2)} \quad&\text{if}\quad T+T^\varepsilon<-t\le2T,
\\
 (-t)^{-\delta(\sigma-1/2)} \quad&\text{if}\quad -t>2T.
\end{cases}
\end{align*}
 By this and the expression \eqref{termC} we get
\begin{align}\label{c-t}
&C''-\frac{T}{2\pi}\zeta(2\sigma)\frac{\zeta'}{\zeta}(\sigma+1/2-it)
\\&\ll
\begin{cases}
T^{2\delta}+ T^{1-\varepsilon\delta(\sigma-1/2)} \quad&\text{if}\quad T+T^\varepsilon<-t\le2T,
\\
(-t)^{2\delta}+T(-t)^{-\delta(\sigma-1/2)} \quad&\text{if}\quad -t>2T.
\end{cases}
 \nonumber
\end{align}

We turn to the term $D''$ defined by the formula \eqref{sumsplit1}.  In view of \eqref{minT} we obtain 
\begin{align}\label{lastProp8}
&D''=-\sum_{(-t-T)^\delta<m\le(-t)^\delta}\sum_{n<m}\frac{-t-m^{1/\delta}}{2\pi}\frac{(m/n)^{it}\Lambda(m/n)}{(mn)^\sigma\sqrt{m/n}}+R_{D''},
\end{align}
where
\begin{align*}
R_{D''}\ll
\begin{cases}
 T^{2\delta} \quad&\text{if}\quad T+T^\varepsilon<-t\le2T,
\\
(-t)^{2\delta} \quad&\text{if}\quad -t>2T.
\end{cases}
\end{align*}
 In formula \eqref{lastProp8} we have that $-t-m^{1/\delta}<T$. Thus
\begin{align}\label{d1d2}
&D''\ll \sum_{(-t-T)^\delta<m\le(-t)^\delta}\sum_{n<m}\frac{T\Lambda(m/n)}{(mn)^\sigma\sqrt{m/n}}+R_{D''}
\\&
=\sum_{n<(-t)^\delta}\sum_{\max(n,(-t-T)^\delta)<m\le(-t)^\delta}\frac{T\Lambda(m/n)}{(mn)^\sigma\sqrt{m/n}}+R_{D''}\nonumber
\\&
=\sum_{n\le(-t-T)^\delta}\sum_{(-t-T)^\delta<m\le(-t)^\delta}\frac{T\Lambda(m/n)}{(mn)^\sigma\sqrt{m/n}}\nonumber
\\&+\sum_{(-t-T)^\delta<n<(-t)^\delta}\sum_{n<m\le(-t)^\delta}\frac{T\Lambda(m/n)}{(mn)^\sigma\sqrt{m/n}}+R_{D''}\nonumber
\\&
=:D''_1+D''_2+R_{D''}.\nonumber
\end{align}
Reasoning as in formulas \eqref{lambdamn} and \eqref{MainTerm}, we get
\begin{align}\label{d1d2ll}
D''_1, D''_2\ll
\begin{cases}
 T^{1-\varepsilon\delta(\sigma-1/2)} \quad&\text{if}\quad T+T^\varepsilon<-t\le2T,
\\
T(-t)^{-\delta(\sigma-1/2)} \quad&\text{if}\quad -t>2T.
\end{cases}
\end{align} 
Formulas \eqref{lastProp8}, \eqref{d1d2}, and \eqref{d1d2ll} yield the bound
\begin{align}\label{dtt} 
D''\ll
\begin{cases}
T^{2\delta}+ T^{1-\varepsilon\delta(\sigma-1/2)} \quad&\text{if}\quad T+T^\varepsilon<-t\le2T,
\\
(-t)^{2\delta}+T(-t)^{-\delta(\sigma-1/2)} \quad&\text{if}\quad -t>2T.
\end{cases}
\end{align}
Summarizing the results obtained in \eqref{acd},  \eqref{rall}, \eqref{c-t}, and \eqref{dtt} we see that
\begin{align}\label{a''-}
&A''-\zeta(2\sigma)\frac{T}{2\pi}\log\frac{T}{2\pi e}-\zeta(2\sigma)\Re\left(\frac{\zeta'}{\zeta}(\sigma+1/2-it)\right)\frac{T}{\pi}
\\&\ll
\begin{cases}
T^{2\delta}+ T^{1-\varepsilon\delta(\sigma-1/2)} \quad&\text{if}\quad T+T^\varepsilon<-t\le2T,
\\
(-t)^{2\delta}+T(-t)^{-\delta(\sigma-1/2)} \quad&\text{if}\quad -t>2T.
\end{cases}\nonumber
\end{align}
Then the formula \eqref{r''} gives that 
\begin{align}\label{r''last}
R''\ll\begin{cases}
 T^{1-\varepsilon\lambda/2} \quad&\text{if}\quad T+T^\varepsilon<-t\le2T,
\\
T(-t)^{-\lambda/2}+T^{1/2}(-t)^{\delta-\lambda/2} \quad&\text{if}\quad -t>2T.
\end{cases}
\end{align}
In view of formulas \eqref{AProp7}, \eqref{a''-}, and \eqref{r''last} we have
\begin{align}\label{lastsum}
&\sum_{0<\gamma\le T}|\zeta(\sigma+i\gamma+it)|^2
-\zeta(2\sigma)\frac{T}{2\pi}\log\frac{T}{2\pi e}-\zeta(2\sigma)\Re\left(\frac{\zeta'}{\zeta}(s+1/2)\right)\frac{T}{\pi}
\\&\ll
\begin{cases}
 T^{2\delta}+T^{1-\varepsilon\delta(\sigma-1/2)}+T^{1-\varepsilon\lambda/2} \quad\text{if}\quad T+T^\varepsilon<-t\le2T,
\\
(-t)^{2\delta}+T(-t)^{-\delta(\sigma-1/2)}+T(-t)^{-\lambda/2}+T^{1/2}(-t)^{\delta-\lambda/2} \ \text{if}\ -t>2T.
\end{cases}\nonumber
\end{align}
From the last formula, replacing   $\lambda/2$ with $\lambda$,  we obtain Proposition \ref{0<-t<T} for $ -t>2T$. For $0<-t\le 2T$, 
  Proposition \ref{0<-t<T} follows by  formulas \eqref{sumsplit}, \eqref{firstsum}, \eqref{secondsum},  and \eqref{lastsum} choosing appropriate constants $\varepsilon$ and $\delta$. This finishes the proof.

\endproof

\proof[Proof of Theorems \ref{prop1} and \ref{th2}.] The theorems immediately follow by Propositions \ref{Ttpositive} and \ref{0<-t<T}.

\endproof

\section{Concluding remarks}\label{conclrem}

Laaksonen and Petridis \cite{lp} investigated a similar sum to the sum of Theorems \ref{prop1} and \ref{th2}. Let $L(s,\chi_1)$ and $L(s,\chi_2)$ be the Dirichlet $L$-functions attached to the primitive Dirichlet characters $\chi_1$ and $\chi_2$. For some fixed prime $P$, let
 $$B(s, P)=\prod_{p\le P}(1-\chi_1(p)p^{-s})(1-\chi_2(p)p^{-s}),$$
 where the product runs over prime numbers $p$.  Under RH, for fixed $1/2<\sigma<1$, they proved that
  \begin{equation*}
            \sum_{0<\gamma\leq T}B(\sigma+i\gamma, P)L(\sigma+i\gamma,\chi_{1})\overline{L(\sigma+i\gamma,\chi_{2})}
            \sim N(T)\sum_{n=1}^{\infty}\frac{d_{n}\overline{\chi}_{2}(n)}{n^{2\sigma}},
        \end{equation*}
        and
        \begin{equation*}
            \sum_{0<\gamma\leq T}B(\sigma+i\gamma,P)\overline{L(\sigma+i\gamma,\chi_{1})}L(\sigma+i\gamma,\chi_{2})\sim N(T)\sum_{n=1}^{\infty}\frac{e_{n}\overline{\chi}_{1}(n)}{n^{2\sigma}},
        \end{equation*}
        where
        \begin{equation*}
            B(s,P)L(s,\chi_{1})=\sum_{n=1}^{\infty}\frac{d_{n}}{n^{s}},\quad  B(s,P)L(s,\chi_{2})=\sum_{n=1}^{\infty}\frac{e_{n}}{n^{s}}.
        \end{equation*}
From this Laaksonen and Petridis \cite{lp} derived the result that, under RH, for a positive proportion of non-trivial zeros of $\zeta(s)$ with $\gamma>0$, the values of the Dirichlet $L$-functions $L(\sigma+i\gamma,\chi_{1})$ and $L(\sigma+i\gamma,\chi_{2})$ are linearly independent over $\mathbb R$.

The discrete mean value of the Dirichlet $L$-function at nontrivial zeros of another Dirichlet $L$-function were investigated by Garunk\v stis and Kalpokas \cite{gk}.  See also Fujii \cite{fuj1, b13}, Conrey, Ghosh and Gonek \cite{b111, b112}, Steuding \cite{b3}, and Garunk\v stis, Kalpokas, and Steuding \cite{gks}.

\end{document}